\documentclass[letterpaper, 10 pt, conference, sort&compress]{ieeeconf}  

\IEEEoverridecommandlockouts                              

\usepackage{amsmath,amssymb,amsfonts}
\usepackage{graphicx}
\usepackage{cite}

\usepackage[version=4]{mhchem}

\usepackage{mathtools}
\usepackage{blkarray, bigstrut}

\usepackage{epstopdf}

\allowdisplaybreaks

\title{\LARGE \bf
Modelling and Economic Optimal Control \\
for a Laboratory-scale Continuous Stirred Tank Reactor \\
for Single-cell Protein Production
}

\author{Marcus Krogh Nielsen, Jens Dynesen, Jess Dragheim,  Ib Christensen, Sten Bay Jørgensen, \\ Jakob Kjøbsted Huusom, Krist V. Gernaey, John Bagterp J{\o}rgensen
\thanks{$^*$This work is funded by Innovation Fund Denmark (9065-00269B). Marcus Krogh Nielsen and John Bagterp J{\o}rgensen are with the Department of Applied Mathematics and Computer Science and Sten Bay J{\o}rgensen, Jakob Kj{\o}bsted Huusom, and Krist V. Gernaey are with the Department of Chemical and Biochemical Engineering at the Technical University of Denmark, Kgs. Lyngby, Denmark. Jens Dynesen, Jess Dragheim, Ib Christensen, and Marcus Krogh Nielsen are with Unibio A/S, Roskilde, Denmark. \texttt{email: \{mkrni,jbjo\}@dtu.dk}}
\thanks{}
}

\begin{document}

\maketitle
\thispagestyle{empty}
\pagestyle{empty}

\begin{abstract}
	In this paper, we present a novel kinetic growth model for the micro-organism \textit{Methylococcus capsulatus} (Bath) that couples growth and pH. We apply growth kinetics in a model for single-cell protein production in a laboratory-scale continuous stirred tank reactor inspired by a physical laboratory fermentor. The model contains a set of differential algebraic equations describing growth and pH-dynamics in the system. We present a method of simulation that ensures non-negativity in the state and algebraic variables. Additionally, we introduce linear scaling of the algebraic equations and variables for numerical stability in Newton's method. Finally, we conduct a numerical experiment of economic optimal control for single-cell protein production in the laboratory-scale reactor. The numerical experiment shows non-trivial input profiles for biomass growth and pH tracking.

\end{abstract}

\section{Introduction}
\label{sec:Introduction}

% ...

Single-cell protein (SCP) provides an alternative source of protein for the feed and food industries to meet the growing demand for protein in the coming decades \cite{oecdfao:2022}. Methanotrophs are bacteria capable of metabolising methane as their source of carbon. Methane is a cheap source of carbon. \textit{Methylococcus capsulatus} (Bath) are methanotrophic bacteria with high protein content that is well-suited for production of SCP \cite{villadsen:etal:2011}.

\textit{M. capsulatus} is aerobic and growth therefore involves fixation of methane and oxygen gas. A U-loop bioreactor has been developed for SCP production. It has been demonstrated that the U-loop bioreactor have good mixing and mass transfer properties \cite{petersen:etal:2017}. Early work on kinetic modelling for growth of \textit{M. capsulatus} shows process instability for high biomass concentrations \cite{olsen:etal:2010a, olsen:etal:2010b, drejer:etal:2017}. Advanced process control techniques have been applied in numerical experiments for monitoring \cite{ritschel:etal:2019a} and to increase stability and productivity in both open- and closed-loop \cite{ritschel:etal:2019b,ritschel:etal:2020}. The metabolism of \textit{M. capsulatus} is well-described in the literature and a metabolic network model exists \cite{lieven:etal:2018}. Metabolic knowledge and process data show that the micro-organism can metabolise several nitrogen sources: ammonium, nitrate, nitrite, and molecular nitrogen. A kinetic growth model exists for growth on these nitrogen sources with methane as carbon source \cite{petersen:etal:2019}. Kinetic parameters for the model are estimated from experimental data \cite{petersen:etal:2020}. The micro-organism metabolises ammonium rather than ammonia directly, pH-dynamics are therefore highly relevant to model precise concentrations of growth-related chemical components. A model including equilibrium reaction in the growth of \text{M. capsulatus} exists, but the kinetics do not depend explicitly on the algebraic variables \cite{prado:etal:2010}. To the best of our knowledge, models coupling growth kinetics and pH-dynamics for \textit{M. capsulatus} as presented in this work have not yet been introduced in the literature.

\begin{figure}[!t]
    \centering
    \includegraphics[width=0.4\textwidth,trim=100 450 100 650,clip]{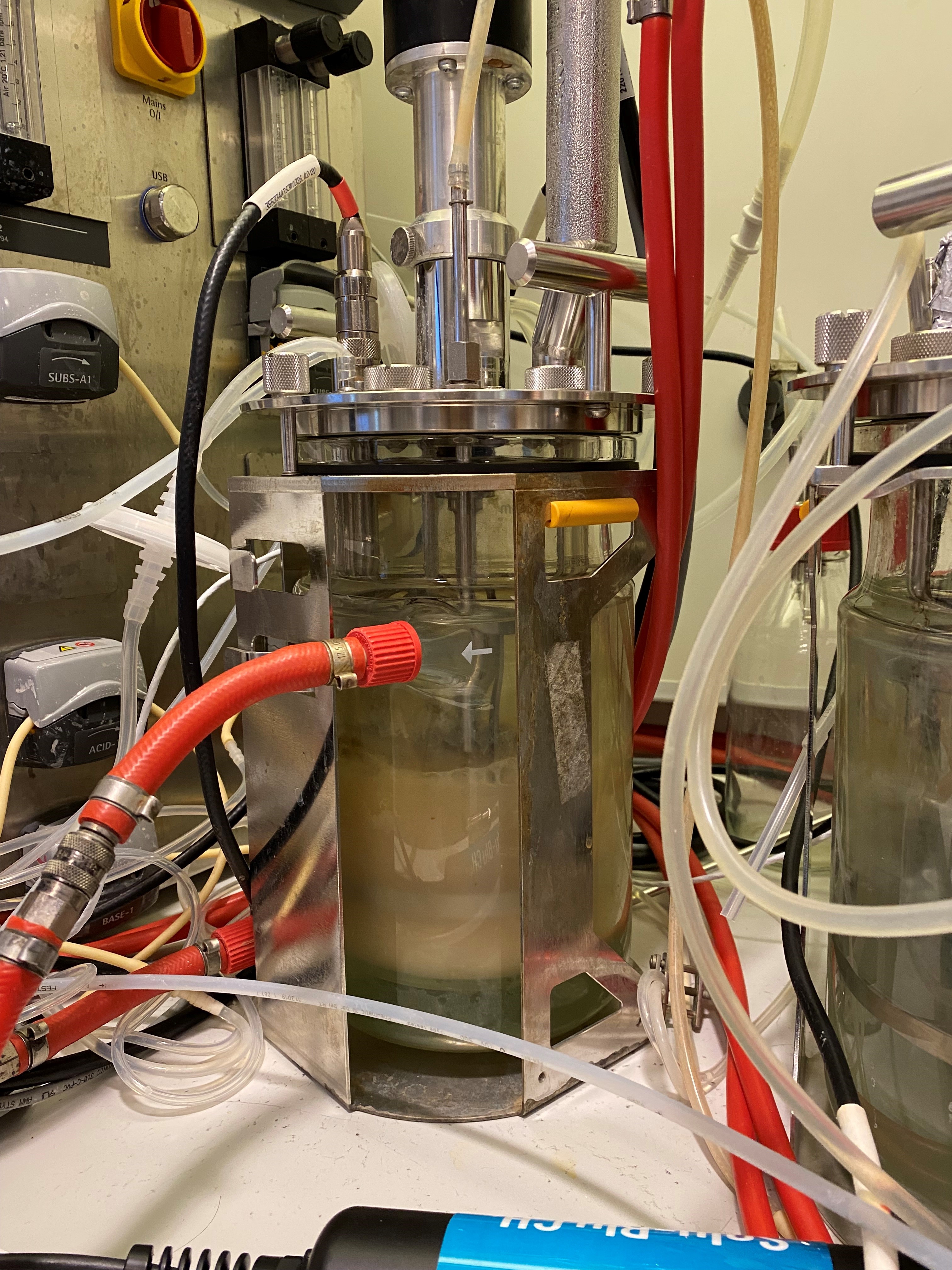}
    \caption{Laboratory-scale reactor for fermention of \textit{M. capsulatus}.}
    \label{fig:cstr}
\end{figure}

In this paper, we present a novel growth model for SCP production by cultivation of the micro-organism \textit{M. capsulatus}. In the model, methane is the carbon source and ammonium is the sole nitrogen source, i.e. growth on nitrite, nitrate, and molecular nitrogen is not included in the model. As such, we apply chemical equilibrium reactions to describe relevant pH-dynamics. We model the fermentation process in a laboratory-scale continuous stirred tank reactor (CSTR), based on a physical laboratory-scale fermentor. Fig \ref{fig:cstr} shows the laboratory reactor. The resulting set of differential algebraic equations (DAEs) model reactor dynamics, growth kinetics, and pH equilibrium dynamics. Optimal control, optimal design of experiments, and parameter estimation motivate the development of the model presented in this work.

The paper is organised as follows. Section \ref{sec:modelling} describes modelling of growth kinetics, pH equilibrium dynamics, as well as reactor dynamics for a CSTR. Section \ref{sec:simulation} presents methods for numerical solution of DAEs with non-negativity constraints as well as scaling of algebraic equations and variables for numerical stability. Section \ref{sec:eocp} describes an economic optimal control problem (OCP) for SCP production in a CSTR. Section \ref{sec:example} presents a numerical experiment of optimal biomass production in a laboratory-scale CSTR. Section \ref{sec:conclusion} presents conclusions.

\section{Modelling}
\label{sec:modelling}
In this section, we describe a mathematical model for SCP production in a CSTR. The model is a system of DAEs in the form
\begin{subequations} \label{eq:model}
    \begin{align}
        \frac{d x}{d t}(t)
            &= f(x(t), y(t), u(t), \theta)
                ,
        &&x(t_{0})
            = x_{0}
                ,   \\
        0
            &= g(x(t), y(t), \theta)
                ,
    \end{align}
\end{subequations}
where $f(\cdot)$ are state dynamical equations, $g(\cdot)$ are algebraic expressions, $t$ is time, $x(t) \in \mathbb{R}^{n_{x}}$ are states, $y(t) \in \mathbb{R}^{n_{y}}$ are algebraic variables, $u(t) \in \mathbb{R}^{n_{u}}$ are manipulated variables, and $\theta$ are system parameters.

\subsection{Variable definitions}
We define the sets
\begin{subequations}
    \begin{align}
        \mathcal{C}_{x,l}
            &= \{ X, N, NO, Na \}
                ,   \\
        \mathcal{C}_{x,g}
            &= \{ S, O, C \}
                ,   \\
        \mathcal{C}_{x}
            &= \mathcal{C}_{x,l} \cup \mathcal{C}_{x,g}
                ,   \\
        \begin{split}
            \mathcal{C}_{y}
                &= \{ \ce{H_{3}O^{+}}, \ce{OH^{-}}, \ce{NH_{3}}, \ce{NH_{4}^{+}},         
                    \\
                &\qquad \ce{CO_{2}}, \ce{H_{2}CO_{3}}, \ce{HCO_{3}^{-}}, \ce{CO_{3}^{2-}} \}
                    ,
        \end{split}
                    \\
        \mathcal{C}_{f,l}
            &= \{ W, N, NO, Na \}
                ,   \\
        \mathcal{C}_{f,g}
            &= \{ S, O \}
                ,   \\
        \mathcal{C}_{f}
            &= \mathcal{C}_{f,l} \cup \mathcal{C}_{f,g}
                .
    \end{align}
\end{subequations}
$\mathcal{C}_{x,l}$ is the set of dissolved state components, $\mathcal{C}_{x,g}$ is the set of state components present as both in gas phase and dissolved in water, $\mathcal{C}_{x}$ is the set of all state components, $\mathcal{C}_{y}$ is the set of algebraic components, $\mathcal{C}_{f,l}$ is the set of liquid inlet flows, $\mathcal{C}_{f,g}$ is the set of gaseous inlet flows, and $\mathcal{C}_{f}$ is the set of all inlet flows. The state components are
\begin{subequations}
    \begin{align}
        X
            &= \ce{CH_{1.8}O_{0.5}N_{0.2}}
                ,   \\
        S
            &= \ce{CH_{4}}
                ,   \\
        O
            &= \ce{O_{2}}
                ,   \\
        N
            &= \{ \ce{NH_{3}}, \ce{NH_{4}^{+}} \}
                ,   \\
        C
            &= \{ \ce{CO_{2}}, \ce{H_{2}CO_{3}}, \ce{HCO_{3}^{-}}, \ce{CO_{3}^{2-}} \}
                ,   \\
        NO
            &= \ce{NO_{3}^{-}}
                ,   \\
        Na
            &= \ce{Na^{+}}
                ,   \\
        S_{g}
            &= \ce{CH_{4}(g)}
                ,   \\
        O_{g}
            &= \ce{O_{2}(g)}
                ,   \\
        C_{g}
            &= \ce{CO_{2}(g)}
                .
    \end{align}
\end{subequations}

\subsection{Stoichiometry and kinetics}
We describe a catabolic reaction and the anabolism of ammonium (not including ATP and ADP) in \textit{M. capsulatus} in terms of the reactions
\begin{small}
    \begin{subequations}
        \begin{align}
            \ce{
            CH_{4} + O_{2} + H_{2}O
                &->
            CO_{2} + 4H_{3}O^{+}
            }  
                ,
            &&r_{1}
                ,   \\
            \ce{
            CH_{4} + O_{2} + \frac{2}{10} NH_{4}^{+}
                &-> 
            X + \frac{15}{10} H_{2}O
            }
                ,
            &&r_{2}
                ,
        \end{align}
    \end{subequations}
\end{small} \vspace{-10pt} \\
where $X = \ce{CH_{1.8}O_{0.5}N_{0.2}}$ is the biomass and $r_{i}(c)$ for $i \in \{1, 2\}$ are the reaction rates. The modelled stoichiometry is
\begin{small}
    \begin{subequations} \label{eq:modelledStoichiometry}
        \begin{align} 
            \ce{
            S + O
                &->
            C
            }
                ,
            &&r_{1}(c)
                ,   \\
            \ce{
            S + O + 0.2 N
                &->
            X
            }
                ,
            &&r_{2}(c)
                .
        \end{align}
    \end{subequations}
\end{small}
The reactions in \eqref{eq:modelledStoichiometry} correspond to the stoichiometric matrix
\begin{equation} \label{eq:stoichiometricMatrix}
    S
        =
    \begin{blockarray}{*{10}{c} l}
        \begin{block}{*{10}{>{$\footnotesize}c<{$}} l}
			X & S & O & N & C & NO & Na & S$_{g}$ & O$_{g}$ & C$_{g}$ \\
        \end{block}
        \begin{block}{[*{10}{c}]>{$\footnotesize}l<{$}}
            0 & -1 & -1 &  0   &  1 &  0 &  0 &  0 &  0 &  0
                \\
            1 & -1 & -1 & -0.2 &  0 &  0 &  0 &  0 &  0 &  0
            	\\
        \end{block}
    \end{blockarray}
\end{equation}
The reaction rates are
\begin{align}
    r_{i}(c)
        &= \mu_{i}(c) c_{X}
            ,
    &&i
        \in \{ 1, 2 \}
            .
\end{align}
The growth rates are
\begin{subequations}
    \begin{align}
        \mu_{1}(c)
            &= \left( \frac{\alpha}{2 \delta} + \frac{8}{20} \right) \mu_{2}(c) + \frac{m}{2 \delta}
                ,   \\
        \mu_{2}(c)
            &= \mu_{\max} \mu_{S}(c) \mu_{O}(c) \mu_{N}(c)
                ,
    \end{align}
\end{subequations}
where $\mu_{\max}$ is the maximum growth rate, $\alpha$ and $\delta$ are parameters related to ATP production and consumption, and $m$ is a maintenance constant. The specific growth rates for methane, oxygen, and ammonium are
\begin{subequations}
    \begin{align}
        \mu_{S}(c)
            &= \frac{c_{S}}{K_{S} \left( 1 + \frac{c_{NH_{4}^{+}}}{K_{N,ox}} \right) + c_{S}}
                ,   \\
        \mu_{O}(c)
            &= \frac{c_{O}}{K_{O} + c_{O}}
                ,   \\
        \mu_{N}(c)
            &= \frac{c_{NH_{4}^{+}}}{K_{N} + c_{NH_{4}^{+}}}
                ,
    \end{align} 
\end{subequations}
The concentrations of dissolved methane, oxygen, and ammonium are denoted $c_{i}$ for $i \in \{S, O, NH_{4}^{+}\}$, respectively. $K_{i}$ for $i \in \{S, O, N\}$ are kinetic constants and $K_{N,ox}$ is the ammonium inhibition constant. We note here that the specific growth for methane, $\mu_{S}(\cdot)$, directly depends on the concentration of an algebraic component, $c_{NH_{4}^{+}}$.

\subsection{pH-balance}
A set of chemical equilibrium reactions, stoichiometric reactions, and mass and charge balances describe the pH equilibrium dynamics in the system.

\paragraph*{Weak acids and bases}
The reactions describing weak acid and base equilibrium are
\begin{small}
    \begin{subequations}
            \begin{align}
                \begin{split}
                   \ce{
                    &2H_{2}O(aq)
                        <=>    \\
                    &\hspace{30pt} H_{3}O^{+}(aq) + OH^{-}(aq)
                    }   
                            ,
                \end{split}
                &&K_{e,W}
                        ,   \\
                \begin{split}
                    \ce{
                    &NH_{3}(aq) + H_{2}O(aq)
                        <=>    \\
                    &\hspace{30pt} NH_{4}^{+}(aq) + OH^{-}(aq)
                    }   
                            ,
                \end{split}
                &&K_{e,N}
                        ,   \\
                \begin{split}
                    \ce{
                    &CO_{2}(aq) + H_{2}O(aq)
                        <=>    \\
                    &\hspace{30pt} H_{2}CO_{3}(aq)
                    }       
                        ,
                \end{split}
                &&K_{e,C1}
                        ,   \\
                \begin{split}
                    \ce{
                    &H_{2}CO_{3}(aq) + H_{2}O(aq)
                        <=>    \\
                    &\hspace{30pt} HCO_{3}^{-}(aq) + H_{3}O^{+}(aq)
                    }       
                        ,
                \end{split}
                &&K_{e,C2}
                        ,   \\
                \begin{split}
                    \ce{
                    &HCO_{3}^{-}(aq) + H_{2}O(aq)
                        <=>    \\
                    &\hspace{30pt} CO_{3}^{2-}(aq) + H_{3}O^{+}(aq)
                    }      
                        ,
                \end{split}
                &&K_{e,C3}
                        .
            \end{align}
    \end{subequations}
\end{small}
The equilibrium constants for water, $K_{e,W}$, ammonia, $K_{e,N}$, and carbon dioxide, $K_{e,C1}$, $K_{e,C2}$, and $K_{e,C3}$, govern the equilibrium reactions.

\paragraph*{Strong acids and bases}
The reactions describing strong acid and base dynamics are
\begin{small}
    \begin{subequations}
        \begin{align}
             \ce{
            NaOH(aq)
                &->
            Na^{+}(aq) + OH^{-}(aq)
            }       
                    ,   \\
            \ce{
            HNO_{3}(aq) + H_{2}O(aq)
                &-> 
            NO_{3}^{-}(aq) + H_{3}O^{+}(aq)
            }       
                    .
        \end{align}
    \end{subequations}
\end{small}
The reactions are assumed to be instantaneous and complete.

\paragraph*{System of algebraic equations}
The chemical equilibrium reactions give rise to the algebraic equations
\begin{subequations}
    \begin{align}
        0
            &= \ce{[H_{3}O^{+}]} \ce{[OH^{-}]} - K_{W}
                ,   \\
        0
            &= \ce{[NH_{4}^{+}]} \ce{[OH^{-}]} - K_{NH} \ce{[NH_{3}]}
                ,   \\
        0
            &= \ce{[H_{2}CO_{3}]} - K_{C,1} \ce{[CO_{2}]}
                ,   \\
        0
            &= \ce{[HCO_{3}^{-}]} \ce{[H_{3}O^{+}]} - K_{C,2} \ce{[H_{2}CO_{3}]}
                ,   \\
        0
            &= \ce{[CO_{3}^{2-}]} \ce{[H_{3}O^{+}]} - K_{C,3} \ce{[HCO_{3}^{-}]}
                .
    \end{align} 
\end{subequations}
The mass balance equations are
\begin{subequations}
    \begin{align}
        0
            &= \ce{[NH_{3}]} + \ce{[NH_{4}^{+}]} - c_{N}
                ,   \\
        \begin{split}
            0
                &= \ce{[CO_{2}]} + \ce{[H_{2}CO_{3}]} + \ce{[HCO_{3}^{-}]} 
                        \\
                &\quad + \ce{[CO_{3}^{2-}]} + \ce{[CO_{2}(g)]} - c_{C}
                    .
        \end{split}
    \end{align}
\end{subequations}
The charge balance equations are
\begin{align}
    \begin{split}
        0
            &= \ce{[OH^{-}]} + \ce{[HCO_{3}^{-}]} - 2 \ce{[CO_{3}]} 
                    \\
            &\quad - c_{NO} - \ce{[H_{3}O^{+}]} - \ce{[NH_{4}^{+}]} - c_{Na}
                .
    \end{split}
\end{align}
We denote the resulting algebraic system of 8 equations and 8 variables
\begin{align}
    0
        &= g(x(t), y(t), \theta)
            .
\end{align}

\subsection{Continuous stirred tank reactor}
We model the laboratory-scale CSTR in the general form
\begin{align}
    \frac{d c_{x}}{d t}
        &= \frac{1}{V} \left( C_{In} - I_{l} c_{x} e_{l}^{T} - I_{g} c_{x} e_{g}^{T} \right) F + Q(c)
            ,
\end{align}
where $c_{x} \in \mathbb{R}^{n_{x}}$ are concentration of the state components, i.e. $c_{i}$ for $i \in \mathcal{C}_{x}$, $V$ is the constant reactor volume, $C_{In} \in \mathbb{R}^{n_{x} \times n_{u}}$ are inlet concentrations of each inlet flow, $F \in \mathbb{R}^{n_{u}}$ are inlet flows, and $c$ are concentrations of all state and algebraic components in the system. The vectors $e_{l} \in \{ 0, 1 \}^{n_{u,l}}$ and $e_{g} \in \{ 0, 1 \}^{n_{u,g}}$ are indicator vectors for liquid and gaseous inflows respectively, i.e. $e_{l}$ is $1$ for all liquid flows and $0$ for all gaseous flows and $e_g$ is $0$ for all liquid flows and $1$ for all gaseous flows. Similarly, $I_{l} = \text{diag}(e_{l})$ and $I_{g} = \text{diag}(e_{g})$ are indicator matrices. The production and mass transfer is defined for components in aqueous solution for $i \in \mathcal{C}_{x,l}$ as
\begin{align}
    Q_{i}(c)
        &= R_{i}(c)
            ,
\end{align}
and for $i \in \mathcal{C}_{x,g}$ as
\begin{align}
    Q_{i}(c)
        &= R_{i}(c) + \frac{1}{1 - \epsilon} J_{gl,i}(c)
            ,
\end{align}
and for components in gas-phase for $i \in \mathcal{C}_{x,g}$ as
\begin{align}
    Q_{i}(c)
        &= - \frac{1}{\epsilon} J_{gl,i}(c)
            .
\end{align}
The production rate is
\begin{align}
    R(c)
        &= S^{T} r(c)
            ,
\end{align}
where $S \in \mathbb{R}^{n_{r} \times n_{x}}$ is the stoichiometric matrix \eqref{eq:stoichiometricMatrix}. The gas-liquid mass transfers for $i \in \mathcal{C}_{x,g}$ are defined as
\begin{align}
    J_{gl,i}(c)
        &= k_{L}a_{i} \left( c_{Sat,i} - c_{i} \right)
            ,
    &&c_{Sat,i}
        = \gamma_{i} c_{i,g}
            ,
\end{align}
The gas-liquid volume fraction is $\epsilon = F_{g}/(F_{l} + F_{g})$, $k_{L}a_{i}$ are mass transfer coefficients, $c_{Sat,i}$ are saturation concentrations, and $\gamma_{i}$ gas-liquid ratio. We apply Henry's law to compute the gas-liquid ratio
\begin{align}
    \gamma_{i}
        &= \frac{R T}{H_{i}^{pc}}
            ,
\end{align}
Henry's constants for chemicals components dissolved in aqueous solution are described in \cite{sander:2015}. Table \ref{tab:modelParameters} describes the model parameters.

\begin{table}[t]
    \centering
    \caption{Model parameters for single-cell protein production in continuous stirred tank reactor.}
    \label{tab:modelParameters}
    \begin{tabular}{l l l | l l l}
        \textbf{name}   & \textbf{value}            & \textbf{unit} 
            &
        \textbf{name}   & \textbf{value}            & \textbf{unit} 
            \\  \hline
        $\mu_{\max}$    & $2.28\mathrm{e}{-}1$      & 1/h
            &
        $m$             & $9.80\mathrm{e}{-}5$      & 1/h
            \\
        $\alpha$        & $2.00\mathrm{e}{-}2$      & -
            &
        $\delta$        & $2.00\mathrm{e}{-}2$      & -
            \\
        $K_{S}$         & $7.50\mathrm{e}{-}5$      & mol/L
            &
        $K_{N,ox}$      & $3.30\mathrm{e}{-}3$      & mol/L
            \\
        $K_{O}$         & $5.50\mathrm{e}{-}5$      & mol/L
            &
        $K_{N}$         & $1.30\mathrm{e}{-}3$      & mol/L
            \\
        $K_{e,W}$       & $1.00\mathrm{e}{-}14$     & -
            &
        $K_{e,N}$       & $5.62\mathrm{e}{-}10$     & -
            \\
        $K_{e,C1}$      & $1.58\mathrm{e}{-}7$      & -
            &
        $K_{e,C2}$      & $4.27\mathrm{e}{-}7$      & -
            \\
        $K_{e,C3}$      & $4.79\mathrm{e}{-}11$     & -
            &
        $c_{In,N}$      & $5.88$                    & mol/L
            \\
        $c_{In,Na}$     & $1.00$                    & mol/L
            &
        $c_{In,NO}$     & $1.00$                    & mol/L
            \\
        $c_{S,g}$       & $1.90\mathrm{e}{-}1$      & mol/L
            &
        $c_{O,g}$       & $1.90\mathrm{e}{-}1$      & mol/L
            \\
        $k_{L}a_{S}$    & $3.89\mathrm{e}{+}2$      & 1/h
            &
        $H^{pc}_{S}$    & $7.05\mathrm{e}{+}2$      & atm/M
            \\
        $k_{L}a_{O}$    & $3.71\mathrm{e}{+}2$      & 1/h
            &
        $H^{pc}_{O}$    & $7.59\mathrm{e}{+}2$      & atm/M
            \\
        $k_{L}a_{C}$    & $3.26\mathrm{e}{+}2$      & 1/h
            &
        $H^{pc}_{C}$    & $2.99\mathrm{e}{+}1$      & atm/M
            \\
        $R$             & $8.21\mathrm{e}{-}2$      & atm/(K M)
            &
        $T$             & $3.15\mathrm{e}{+}2$      & K
            \\
        $V$             & $1.00$                    & L
            &
        -               & -                         & -
    \end{tabular}
\end{table}

\section{Simulation}
\label{sec:simulation}
In this section, we describe Euler's implicit method for numerical solution of DAEs. Additionally, we describe a variation Newton's method, which yields only non-negative solutions.

\subsection{Euler's implicit method}
For differential algebraic systems in the form presented in \eqref{eq:model}, we may discretise with Euler's implicit method as
\begin{subequations} \label{eq:implicitEuler}
    \begin{align} 
        x_{k+1}
            &= x_{k} + f(x_{k+1}, y_{k+1}, u_{k}, \theta) \Delta t
                ,   \\
        0
            &= g(x_{k+1}, y_{k+1}, \theta)
                .
    \end{align}
\end{subequations}
From this, we define the residual function
\begin{align} \label{eq:implicitEulerResidual}
    \begin{aligned} 
        R(z_{k+1})
            &= R(x_{k+1}, y_{k+1})
                    \\
            &=
        \begin{bmatrix}
            x_{k+1} - x_{k} - f_{k}(x_{k+1}, y_{k+1}) \Delta t
                \\
            g_{k}(x_{k+1}, y_{k+1})
        \end{bmatrix}
                ,
    \end{aligned}
\end{align}
where $z_{k} = [x_{k}; y_{k}]$, $f_{k}(x, y) = f(x, y, u_{k}, \theta)$ and $g_{k}(x, y) = g(x, y, \theta)$. Solutions to the difference equation, \eqref{eq:implicitEuler}, are roots in the residual equation, \eqref{eq:implicitEulerResidual}, obtained by solving
\begin{align}
    0
        &= R(z_{k+1})
            .
\end{align}
For chemical systems modelling concentrations, only non-negative solutions represent solutions in the physical system.

\subsection{Newton's method with asbolute step}
Newton's method provides a means of iteratively finding roots given a function and its Jacobian. For a residual function, $R(z)$, the system of equations
\begin{align}
    0
        &= R(z)
            ,
\end{align}
can be solved iteratively, where the search direction, $\Delta z$, is obtained as the solution of
\begin{align}
    0
        &= \frac{\partial R}{\partial z}(z) \Delta z + R(z)
            .
\end{align}
The Newton step in each iteration, $n$, is
\begin{align}
    z_{n+1}
        &= z_{n} + \Delta z_{n}
            ,
    &&\Delta z_{n}
        = - \left( \frac{\partial R}{\partial z}(z_{n}) \right)^{-1} R(z_{n})
            .
\end{align}
Implicit Euler discretisation of DAE systems in the form \eqref{eq:model} results in the residual Jacobian,
\begin{align} \label{eq:jacobian}
    \begin{aligned}
        \frac{\partial R}{\partial z}(z)
            &= 
            \begin{bmatrix}
                I - \frac{\partial f_{k}}{\partial x} (x, y) \Delta t     & - \frac{\partial f_{k}}{\partial y}(x, y) \Delta t
                    \\
                \frac{\partial g_{k}}{\partial x}(x, y)                  & \frac{\partial g_{k}}{\partial y}(x, y)
            \end{bmatrix}
                .
    \end{aligned}
\end{align}
In the case where only non-negative solutions to the residual also represent solutions to the physical system, we may apply so-called absolute Newton's method, such that the step becomes
\begin{align}
    z_{n+1}
        &= \left| z_{n} + \Delta z_{n} \right|
            .
\end{align}
This variation of Newton's method ensures non-negativity in the roots, but does not maintain the convergence properties of the original formulation of the method. The method is empirical, but has been successfully applied to chemical equilibrium systems  \cite{meintjes:morgan:1987}.

\subsection{Scaling algebraic equations and variables}
The different scales of the algebraic equations and variables may result in ill-conditioning in the Jacobian of the residual \eqref{eq:jacobian}. This can lead to numerical inaccuracies and divergence in Newton's method. As such, we introduce a scaling of the algebraic variables and equations
\begin{align}
    \tilde{g}(x(t), \tilde{y}(t), \theta)
        &= S_{g} g(x(t), S_{y} \tilde{y}(t), \theta)
            ,   \\
    y(t) 
        &= S_{y} \tilde{y}(t)
            .
\end{align}
The matrices, $S_{g} = \text{diag}(s_{g})$ and $S_{y} = \text{diag}(s_{y})$, are diagonal scaling matrices. The vectors, $s_{g} \in \mathbb{R}^{n_{y}}$ and $s_{y} \in \mathbb{R}^{n_{y}}$, define the scaling factors for the algebraic equations and variables, respectively. The resulting scaled Jacobian is
\begin{align} \label{eq:scaledJacobian}
    \frac{\partial \tilde{g}}{\partial \tilde{y}}
        &= \frac{\partial \tilde{g}}{\partial g} \frac{\partial g}{\partial y} \frac{\partial y}{\partial \tilde{y}}
        = S_{g} S_{y} \frac{\partial g}{\partial y}
            .
\end{align}
We choose $S_{g}$ and $S_{y}$ such that the Jacobian \eqref{eq:scaledJacobian} is well-conditioned.

\section{Economic optimal control problem}
\label{sec:eocp}
In this section, we present a formulation of an economic OCP for SCP production in a CSTR. 

\subsection{Formulation}
Consider the economic OCP
\begin{subequations}
    \begin{align}
        \min_{x(t), y(t), u(t)} \quad
            &\phi = \alpha_{eco} \phi_{eco} + \alpha_{pH} \phi_{pH} + \alpha_{du} \phi_{du}
                ,   \\
        s.t. \quad
            &x(t_{0})
                = x_{0}
                    ,   \\
            &\frac{d x}{d t}(t)
                = f(x(t), y(t), u(t), d(t), \theta)
                    ,   \\
            &0
                = g(x(t), y(t), \theta)
                    ,   \\
            &x_{\min}
                \le x(t) \le x_{\max}
                    ,   \\
            &y_{\min}
                \le y(t) \le y_{\max}
                    ,   \\
            &u_{\min}
                \le u(t) \le u_{\max}
                    .
    \end{align}
\end{subequations}
The economic objective is
\begin{align}
    \phi_{eco}
        &= \phi_{cost} - \phi_{profit} - \phi_{ctg}
            ,
\end{align}
where
\begin{subequations}
    \begin{align}
        \phi_{profit}
            &= \int_{t_{0}}^{T} p_{X} c_{X}(t) e^{T} F_{l}(t) dt
                ,   \\
        \phi_{cost}
            &= \int_{t_{0}}^{T} p_{F}^{T} F(t) dt
                ,   \\
        \phi_{ctg}
            &= p_{X} V \left( c_{X}(T) - c_{X}(t_{0}) \right)
                .
    \end{align}
\end{subequations}
The objective term $\phi_{profit}$ is the value of the harvested biomass over the horizon, $\phi_{cost}$ is the cost associated with the inlet streams of gaseous and liquid substrates, and $\phi_{ctg}$ is the cost-to-go. The cost-to-go term is included such that the reactor is not emptied for additional profit at the end of the control horizon. $p_{X}$ [USD/g] and $p_{F}$ [USD/L] are the unit values of the biomass and inlets, respectively. $F_{l}(t)$ [L/h] are the liquid inlet streams, $e^{T} F_{l}(t)$ [L/h] is the liquid outlet stream, and $V$ [L] is the reactor volume. We define the pH tracking objective term, as
\begin{align}
    \phi_{pH}
        &= \frac{1}{2} \int_{t_{0}}^{T} \left\| \overline{pH} - pH(t) \right\|_{Q_{pH}}^{2} dt
            ,
\end{align}
where $\overline{pH}$ is the target, $pH(t) = - \log_{10} \left( \ce{[H_{3}O^{+}]} \right)$, and $Q_{pH}$ is a weight matrix. We define the input rate-of-movement objective as
\begin{align}
    \phi_{du}
        &= \frac{1}{2} \int_{t_{0}}^{T} \left\| \frac{du}{dt}(t) \right\|_{Q_{du}}^{2} dt
            ,
\end{align}
such that the changes in the manipulated inputs over time, $u(t)$, are penalised quadratically.

\section{Numerical experiment}
\label{sec:example}
\begin{figure}
    \centering
    \includegraphics[width=0.45\textwidth]{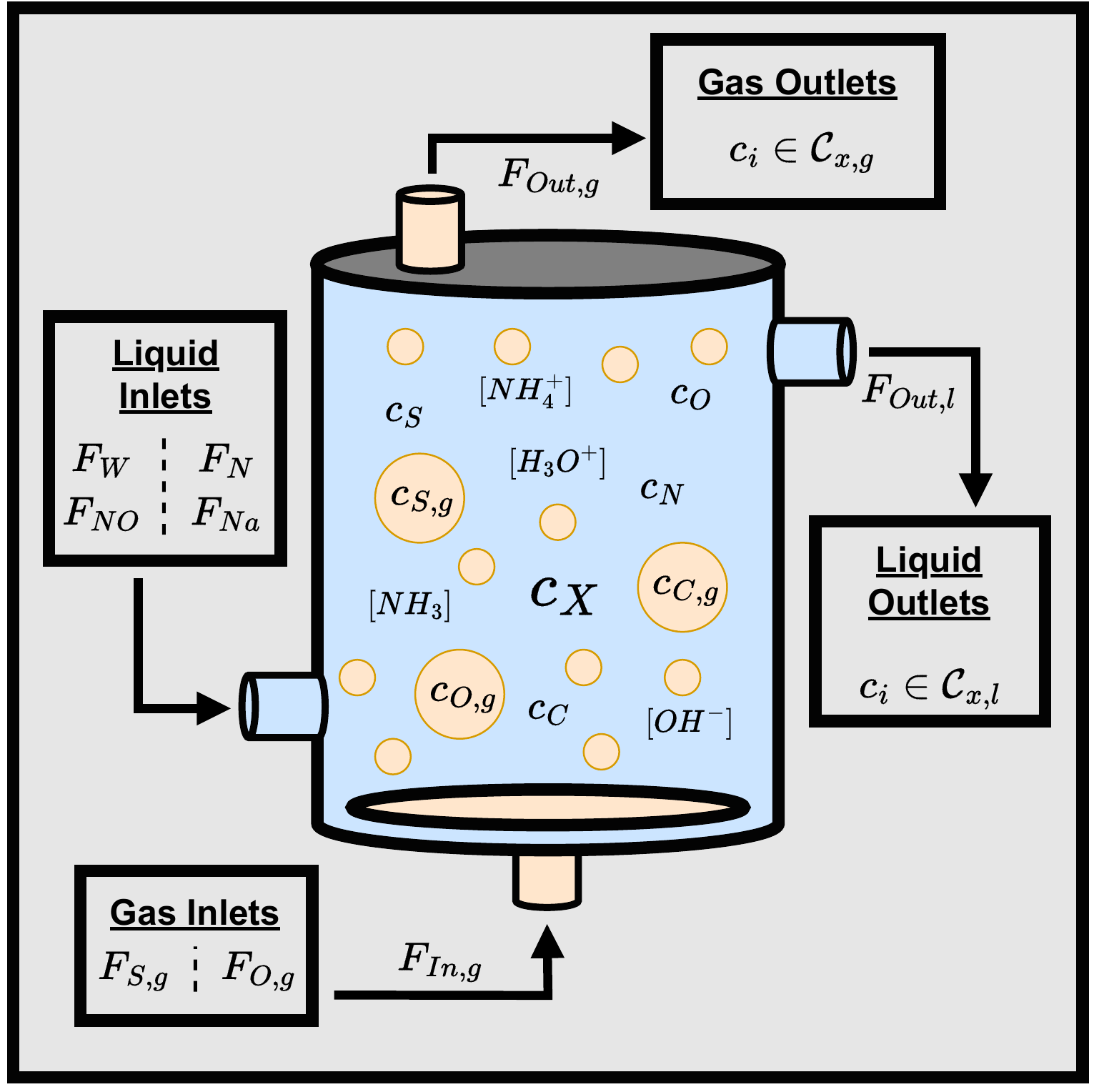}
    \caption{Illustration of the laboratory-scale continuous stirred tank reactor.}
    \label{fig:cstrDiagram}
\end{figure}

In this section, we present a simulation example for economic optimal control in a laboratory-scale CSTR for SCP production. Fig. \ref{fig:cstrDiagram} is an illustration of the fermentor.

\subsection{Laboratory fermentor}
In the physical laboratory fermentor, operation is separated into a batch phase and a continuous phase. 

\paragraph*{Batch phase} 
The reactor is initialised from low biomass concentration, $c_{X} < 0.1$ g/L (inoculum), with substrates present for growth to a biomass concentration of $c_{X} \ge 2$ g/L. In this phase, only gasses flow in and out of the reactor, i.e. $F_{i} = 0$ for $i \in \mathcal{C}_{f,l}$ and $F_{j} \ge 0$ for $j \in \mathcal{C}_{f,g}$.

\paragraph*{Continuous phase}
The reactor is initialised from a biomass concentration high enough to consistently sustain continuous operation, i.e. $c_{X} \ge 2$ g/L. In the continuous phase, liquid substrates flow into the reactor. To keep the volume constant, liquid reactor content is harvested at the same rate as the inflow, i.e. $F_{i} \ge 0$ for $\mathcal{C}_{f,g}$, $F_{j} \ge 0$ for $i \in \mathcal{C}_{f,l}$, and $F_{Out,l} = \sum_{j \in \mathcal{C}_{f,l}} F_{j}$.

\subsection{Implementation}
We formulate the economic OCP for the continuous phase of the laboratory-scale CSTR. We implement the discretised OCP in Matlab with a simultaneous collocation-based approach. We apply the implicit Euler discretisation scheme (i.e. right rectangular rule). Casadi is applied to solve the resulting nonlinear program \cite{anderson:etal:2019}. The initial state is $x_{0} = [ 2.00, 2.31\mathrm{e}{-}2, 3.77\mathrm{e}{-}2, 4.03\mathrm{e}{-}1, 9.10\mathrm{e}{-}1, 3.07\mathrm{e}{-}3, \\ 2.23\mathrm{e}{-}7, 6.29\mathrm{e}{-}1, 1.11, 1.05 ]^{T}$ [g/L] and the initial input is $u_{0} = [ 1.00\mathrm{e}{-}2, 6.84\mathrm{e}{-}5, 6.71\mathrm{e}{-}7, 5.62\mathrm{e}{-}11, 8.96\mathrm{e}{-}3, \\8.53\mathrm{e}{-}3 ]^{T}$ [L/h]. The unit prices for substrates are $p_{W} = 0.00$, $p_{N} = 1.00\mathrm{e}{-}3$, $p_{NO} = 1.00\mathrm{e}{-}1$, $p_{Na} = 1.00\mathrm{e}{-}1$, $p_{S,g} = 1.00\mathrm{e}{-}3$, $p_{O,g} = 1.00\mathrm{e}{-}3$ [USD/L] and for biomass $p_{X} = 1.0\mathrm{e}{-}2$ [USD/g]. $\alpha_{eco} = 1.0$, $\alpha_{pH} = 2.0 \mathrm{e}{+}2$, and $\alpha_{du} = 1.0$ are the scaling parameters for the objectives in the OCP. The weight matrices for the pH tracking and input rate-of-movement objectives are $Q_{pH} = I$ and $Q_{du} = \text{diag}([1.0, 1.0, 10.0, 10.0, 0.1, 0.1])$, respectively. The algebraic equations are scaled with $s_{g} = [1.0\mathrm{e}{+}7, 1.0\mathrm{e}{+}7, 1.0\mathrm{e}{+}4, 1.0\mathrm{e}{+}8, 1.0\mathrm{e}{+}10, 1.0\mathrm{e}{+}0, 1.0\mathrm{e}{+}0, \\1.0\mathrm{e}{+}0]$ and variables with $s_{y} = \{ 1 \}^{n_{y}}$. The algebraic variables are parametrised in the implementation, such that $y(\alpha(t)) = 10^{-\alpha(t)}$, i.e. we describe the pH-value directly as $pH(t) = - \log_{10}(y_{H_{3}O^{+}}(t)) = \alpha_{H_{3}O^{+}}(t)$. This parametrisation also ensures non-negative during the optimisation. We solve the discretised OCP over a horizon of 48 hours with $200$ discrete time-steps. All variables have a lower boundary of $0$. The biomass and ammonium concentrations have upper boundaries $c_{X} \le 20.0$ g/L and $c_{N} \le 1.0$ M ($17.03$ g/L), respectively.

\subsection{Results}
\begin{figure*}[t]
    \begin{center}
    \includegraphics[width=0.99\textwidth]{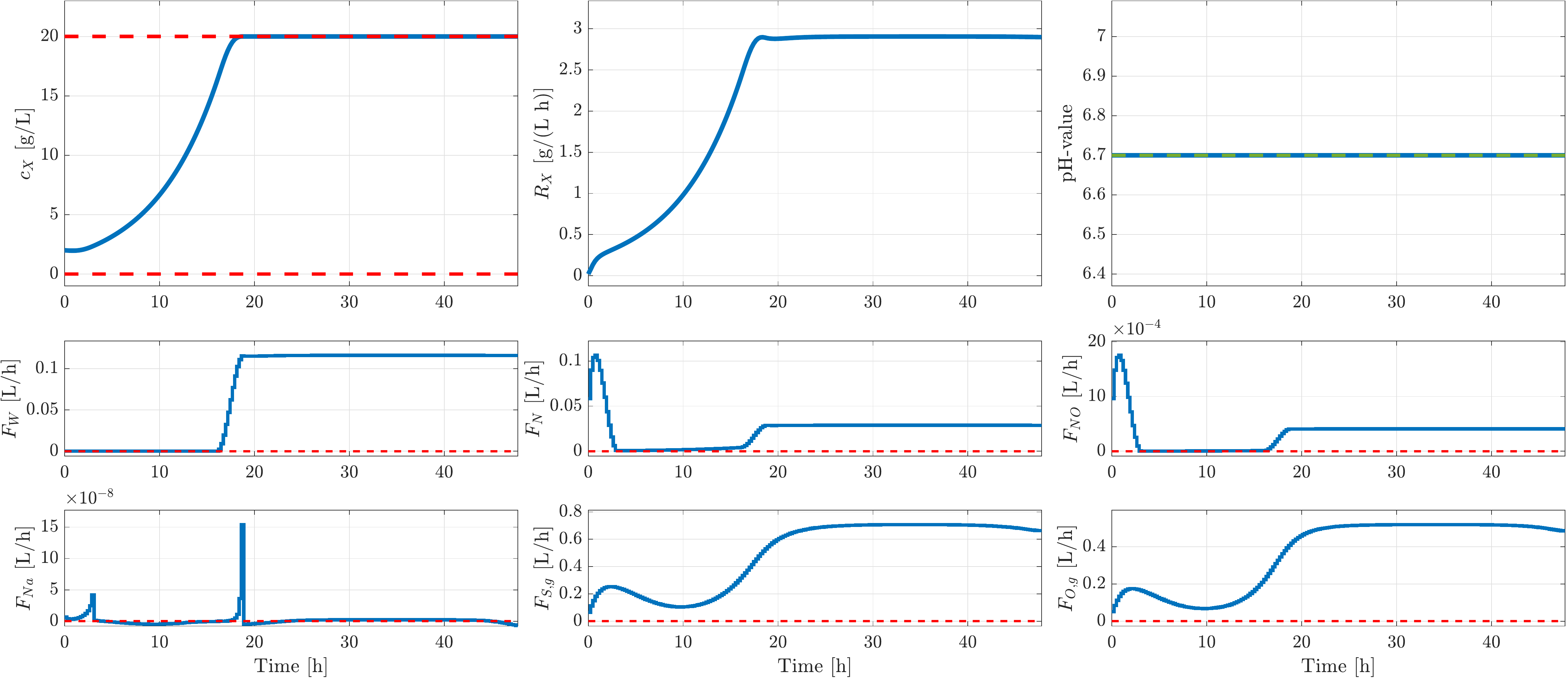}
    \end{center}
    \caption{Solution to economic optimal control problem for single-cell protein production in continuous-stirred tank reactor. The solution shows an initial growth phase of around 18 hours, followed by a production phase with increased inlet flow-rate of water for stable and high productivity.}
    \label{fig:ocpKPI}
\end{figure*}

\begin{figure*}[t]
    \begin{center}
    \includegraphics[width=0.99\textwidth]{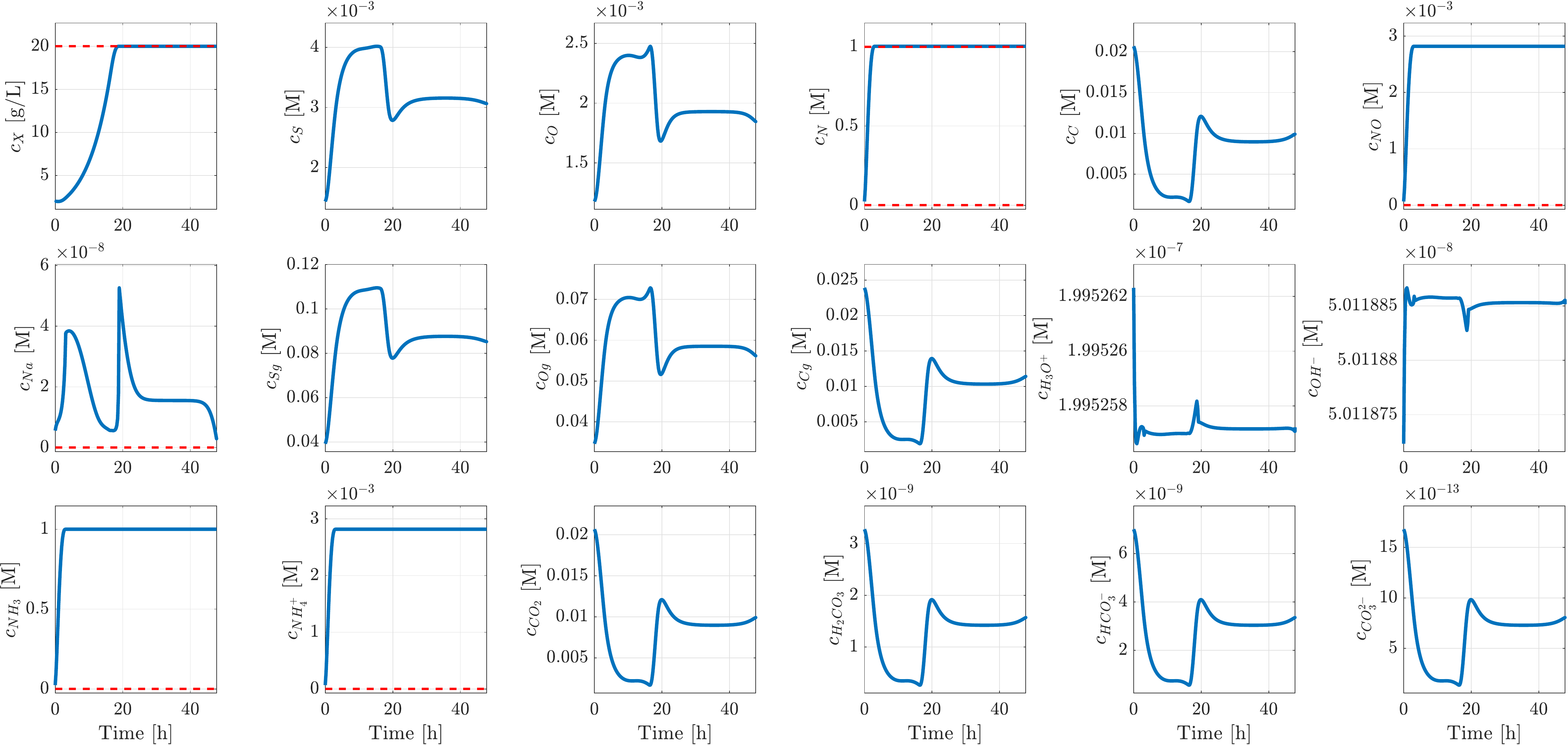}
    \end{center}
    \caption{Concentrations of all state and algebraic variables computed as the solution to the economic optimal control problem.}
    \label{fig:ocpConcentrations}
\end{figure*}

Fig. \ref{fig:ocpKPI} illustrates the biomass concentration, biomass productivity, pH-value, and liquid and gas inlet streams of the numerical experiment. Fig. \ref{fig:ocpConcentrations} illustrates the concentrations of all state and algebraic variables in the numerical experiment. We observe clear separation between growth and production phases in the solution. We note that the nitrogen source, ammonium, is at the upper limit during production. This indicates that the nitrogen source is the limiting substrate during production.

\section{Conclusion}
\label{sec:conclusion}
This paper presented a novel growth model for SCP production in a laboratory-scale CSTR. The DAE model describes reactor dynamics, growth kinetics of the micro-organism \textit{M. capsulatus}, as well as pH equilibrium dynamics in the system. The model couples growth and pH as the growth directly depends on the algebraic variables. Finally, we conducted a numerical experiment of economic optimal control for the system. In the numerical experiment, we demonstrated optimal biomass growth and production, while tacking the pH-value in the reactor.

\bibliographystyle{IEEEtran}
\bibliography{references/bibliography}

\end{document}